\documentclass{amsart}

\usepackage[frame,cmtip,curve,arrow,matrix,line,graph]{xy}
\usepackage{latexsym}
\usepackage{euscript}
\newtheorem{theorem}{Theorem}[section]
\newtheorem{proposition}[theorem]{Proposition}
\newtheorem{corollary}[theorem]{Corollary}

\theoremstyle{definition}
\newtheorem{definition}[theorem]{Definition}
\newtheorem{example}[theorem]{Example}
\newtheorem{remark}[theorem]{Remark}
\newtheorem{remarks}[theorem]{Remarks}

\newcommand{\qq}{\mathbb{Q}}
\newcommand{\cc}{\mathbb{C}}

\newcommand{\zz}{\mathbb{Z}}

\newcommand{\id}{\textrm{id}}
\newcommand{\Hom}[2]{\mathrm{Hom}\left(#1,#2\right)}

\newcommand{\Stab}{\mathrm{Stab} \,}

\newcommand{\git}{/\!\!/}
\newcommand{\udot}{\,{\bf \dot{}}\,}

\newcommand{\sh}[1]{{\bf Sh}(#1)}

\newcommand{\icg}[1]{\mathcal{IC}\udot_{\!\!\!G}(#1)}
\newcommand{\eqc}[2]{\mathcal{C}\udot_{\!\!\!#1}(#2)}
\newcommand{\eqic}[2]{\mathcal{IC}\udot_{\!\!\!#1}(#2)}

\newcommand{\ic}[1]{\mathcal{IC}\udot (#1)}

\newcommand{\dgs}[1]{\mathcal{#1}\udot}
\newcommand{\dgsx}[2]{\mathcal{#1}\udot{(#2)}}
\newcommand{\shf}[1]{\mathcal{#1}}

\newcommand{\rdf}[1]{R#1}

\newcommand{\q}{Q}

\newcommand{\der}[1]{{\bf D}(#1)}
\newcommand{\eqder}[2]{{\bf D}_{#1}(#2)}

\newcommand{\ih}[2]{I\!H^{#1}(#2)}
\newcommand{\eqih}[3]{I\!H^{#1}_{#2}(#3)}
\newcommand{\h}[2]{H^{#1}(#2)}
\newcommand{\eqh}[3]{H^{#1}_{#2}(#3)}

\newcommand{\shom}[2]{H^{#1}(#2)}

\newcommand{\perv}[1]{{\bf Perv}(#1)}
\newcommand{\eqperv}[2]{{\bf Perv}_{#1}(#2)}
\newcommand{\res}[2]{{\bf Res}_{#1}(#2)}
\newcommand{\var}{\bf Var}
\newcommand{\eqsh}[2]{{\bf Sh}_{#1}(#2)}
\newcommand{\forget}{\textrm{For}}

\begin{document}
\author{Jonathan Woolf}
\address{Department of Mathematics and 
Statistics, University of Edinburgh, Edinburgh, EH9~3JZ, Scotland;
jon@maths.ed.ac.uk}
\title{The decomposition theorem and the intersection cohomology of quotients in algebraic geometry}
\date{September 2000}
\subjclass{14F43,14L24}

\begin{abstract}
Suppose a connected reductive complex algebraic group $G$ acts linearly on a complex projective variety $X$. We prove that if 
$$
1 \to N \to G \to H \to 1
$$ 
is a short exact sequence of connected reductive groups and $X^{ss}$ the open set of semistable points for the action of $N$ on $X$ then $\eqih{*}{H}{X^{ss}\git N}$ is a direct summand of $\eqih{*}{G}{X^{ss}}$. The inclusion is provided by the decomposition theorem and certain resolutions of the action allow us to define projections. 
\end{abstract}

\maketitle

\section*{Introduction}
\label{intro}
Suppose a connected complex reductive group $G$ acts linearly on an ample line bundle over a complex projective variety $X$.  There is a geometric invariant theory quotient $X^{ss} \git G$ of the set $X^{ss}$ of semistable points of $X$ for this linearisation. We assume that the open set $X^s$ of stable points for the linearised action is non-empty in order that $\dim X^{ss} \git G = \dim X - \dim G$. When the variety $X$ is smooth, and every semistable orbit is furthermore stable, the relation between the equivariant cohomology of $X$ and the cohomology of $X^{ss}\git G$ has been intensively studied, see eg. \cite{k1,wi,jk,gk,tw} and many others. Aside from the intrinsic interest of relating an equivariant invariant to one defined on a quotient, the subject has obvious applications in the topological study of various moduli spaces. A key point in this theory is the observation that  $H^*_G(X^{ss}) \cong H^*(X^{ss}\git G)$ under these assumptions. Several papers, for instance \cite{k3,h1,lt} and \cite{kiem}, have studied what happens when we relax the assumption that every semistable orbit be stable. Since the quotient $X^{ss}\git G$ will in general then be singular they all consider its intersection cohomology groups rather than its ordinary cohomology. Principally this is because the former retain for singular projective varieties the structures, collectively known as the K\"ahler package, which hold for the cohomology of a smooth projective variety. So, although they are less tractable in terms of functoriality, in many ways they provide a richer invariant for the study of singular varieties. Another common theme in these papers is the construction of some resolution of the $G$ action on $X$ and the use of this to define a map $\eqh{*}{G}{X^{ss}} \to \ih{*}{X^{ss}\git G}$ which is surjective. In this paper we analyse this approach in general terms and provide a framework in which to place their results.

Let us allow that the complex projective variety $X$ may be singular and consider the quotient morphism $X^{ss} \to X^{ss} \git G$. We would like to relate the intersection cohomology $\ih{*}{X^{ss}\git G}$ of the quotient to the equivariant intersection cohomology $\eqih{*}{G}{X^{ss}}$ of the semistable points. Throughout this paper we will take coefficients in the rationals $\qq$ and thus avoid issues of torsion. In slightly different terms we note that $\eqih{*}{G}{X^{ss}}$ is, almost by definition, the intersection cohomology of the quotient stack $[X^{ss}/G]$ so that we are comparing the quotient stack with the geometric invariant theory quotient. When every semistable orbit is stable \cite[thm. 9.1]{bl} tells us that
$$
\eqih{*}{G}{X^{ss}} \cong \ih{*}{X^{ss}\git G}.
$$
This of course corresponds to the fact that the quotient stack $[X^{ss}/G]$ is represented by the orbifold $X^{ss}/G \cong X^{ss}\git G$.
If we remove this condition the picture is more complicated. Intuitively the equivariant intersection cohomology groups contain more information. This notion is expressed quite simply in our main result which says that $\ih{*}{X^{ss}\git G}$ is always a direct summand of $\eqih{*}{G}{X^{ss}}$. In fact we prove something slightly more general. Suppose $N$ is a connected normal reductive subgroup of $G$, and $H$ the corresponding quotient $G/N$. Let $X^{ss}$ now stand for the set of $N$-semistable points of $X$. Then we show that $\eqih{*}{H}{X^{ss}\git N}$ is a direct summand of $\eqih{*}{G}{X^{ss}}$.

It is no great surprise that the key ingredient of our proof is the decomposition theorem of \cite[\S 6]{bbd}, or rather its equivariant extension which was proved in \cite[\S 5]{bl}. To put this result in context we review, in \S \ref{eq der cat}, the construction of the constructible equivariant derived category from \cite{bl}. Then in \S \ref{pull-back} we show how to define a map
\begin{equation}
\label{first}
\eqih{*}{H}{X^{ss}\git N} \to \eqih{*}{G}{X^{ss}}.
\end{equation}
In the final section we consider certain resolutions of $X^{ss}$, which we term $N$-stable, and show how these can be used to construct projections
$$
\eqih{*}{G}{X^{ss}}\to \eqih{*}{H}{X^{ss}\git N} 
$$
and thence show that (\ref{first}) is an inclusion. 
\subsection*{Acknowledgements}
I am greatly indebted to interesting and illuminating conversations with both Young-Hoon Kiem and Sue Tolman. I would also like to thank Frances Kirwan for her patience in reading preliminary drafts.

\tableofcontents

\section{The equivariant derived category}
\label{eq der cat}
The results of this paper are couched in the language of the equivariant derived category introduced by Bernstein and Lunts in \cite{bl}. We give a very brief review of the structures which we use, but refer the reader to \cite{bl} for the details. Background material on the constructible derived category, $t$-structures, perverse sheaves and intersection cohomology can be found in \cite{gm2} and \cite{bbd}.

Suppose that $X$ is a complex algebraic variety. A `sheaf on $X$' will be a module over the constant sheaf with coefficients in $\qq$, and  will furthermore be constructible i.e. its cohomology sheaves will be locally constant on the strata of some stratification of $X$ by smooth subvarieties. Such sheaves form an Abelian category $\sh{X}$. We write $\der{X}$ for the bounded below derived category of constructible sheaves on $X$. The middle perversity intersection cohomology complex $\ic{X}$ is defined up to quasi-isomorphism as an object of $\der{X}$ obeying axioms set out in \cite{gm2}. We follow their definition except that we shift by the complex dimension $d_X$ of $X$ and so follow the Beilinson--Bernstein--Deligne--Gabber indexing used in \cite{bbd,bl} (see \cite[2.3]{gm2} for a comparison of the various indexing systems). In addition to the standard $t$-structure on $\der{X}$, whose heart is $\sh{X}$, there is a $t$-structure associated to the middle perversity whose heart is the Abelian category $\perv{X}$ of perverse sheaves --- see \cite[\S 5]{bl} and \cite[\S 2]{bbd}. $\ic{X}$ is a simple object in this full subcategory in the sense that it has no subobjects (see \cite[4.3]{bbd}). 

A map $\varphi : X \to X'$ gives rise to functors $\varphi^*, \varphi^! : \der{X'} \to \der{X}$ and $\varphi_*,\varphi_! : \der{X} \to \der{X'}$. Here $\varphi^*$ is the left adjoint of $\varphi_*$ and $\varphi_!$ the left adjoint of $\varphi^!$. There is also a natural tensor product $\otimes$ on $\der{X}$.
\begin{remark}
Note that following \cite{bl} we write $\varphi_*$ and not the more usual $\rdf \varphi_*$ for the right derived functor because we will always work with derived categories and so there is no possibility of confusion with the push-forward of sheaves.
\end{remark}
The intersection cohomology groups of $X$ are defined by
$$
\ih{*}{X} := \shom{*+d_X}{\pi_* \ic{X}}
$$
where $\pi: X \to pt$ is the map to a point, and $d_X$ the complex dimension of $X$. 

Now suppose that a reductive algebraic group $G$ acts algebraically on $X$. If $X$ is a principal $G$-space then the derived category $\der{X/G}$ of the quotient is a good definition of `equivariant derived category'. More generally let $\res{G}{X}$ be the category of $G$-resolutions of $X$ i.e. the category whose objects are equivariant morphisms $Y \to X$ where $Y$ is a principal $G$-space and whose morphisms are commutative diagrams
\[
\xymatrix{
Y \ar[dr] \ar[rr] && Y'  \ar[dl]\\
& X
}
\]
of equivariant morphisms. There is a natural functor $\Phi : \res{G}{X} \to \var$ to the category of varieties given by $Y \mapsto Y/G$. Both $\sh{-}$ and $\der{-}$ can naturally be viewed as fibred categories over $\var$. It is shown in \cite[\S 2]{bl} that the fibres of these over the functor $\Phi$ are respectively the category $\eqsh{G}{X}$ of constructible equivariant sheaves and the constructible bounded below equivariant derived category $\eqder{G}{X}$. More concretely an object $\dgs{A}$ of the latter is given by an object $\dgsx{A}{Y} \in \der{Y/G}$ for each $Y \in \res{G}{X}$ and, functorially, for each morphism $\alpha: Y \to Y'$ a quasi-isomorphism $\dgsx{A}{Y} \cong \alpha^* \dgsx{A}{Y'}$. If $X$ is itself a principal $G$-space then the category $\res{G}{X}$ has a final object, namely $X$ itself, and we see that, as hoped, $\eqder{G}{X} \cong \der{X/G}$.

The equivariant derived category inherits the structure of a triangulated category. The assignment $\dgs{A} \mapsto \dgsx{A}{G\times X}$ defines a forgetful functor $\forget : \eqder{G}{X} \to \der{X}$ and both the usual and perverse $t$-structures on $\der{X}$ can be lifted via the forgetful functor to $t$-structures on $\eqder{G}{X}$. The respective hearts are $\eqsh{G}{X}$ and the equivariant perverse sheaves $\eqperv{G}{X}$. Furthermore the intersection cohomology complex $\ic{X}$ has an equivariant lift to a simple object $\eqic{G}{X} \in \eqperv{G}{X}$ --- see \cite[\S 5]{bl}. 

\begin{remarks}
\begin{enumerate}
\item Note that, owing to the distinction between term-by-term isomorphisms and quasi-isomorphisms of complexes, the equivariant derived category $\eqder{G}{X}$ is \emph{not} generally equivalent to the derived category of $\eqsh{G}{X}$.
\item Generalising the statement that $\eqder{G}{X} \cong \der{X/G}$ for a principal $G$-space $X$, we can interpret the equivariant sheaves and the equivariant derived category as $\sh{[X/G]}$ and $\der{[X/G]}$ where $[X/G]$ is the quotient stack. This follows immediately from the definition of stacks via fibred categories, see for example \cite[2.2]{gom}.
\end{enumerate} 
\end{remarks}
The functors $\varphi_*,\varphi^*,\varphi_!,\varphi^!$ and $\otimes$ extend to the equivariant context (where $\varphi$ is now an equivariant map) simply by defining $(\varphi^*\dgs{A})(Y) = \varphi^*(\dgsx{A}{Y})$ etc. There are also two new functors in the equivariant setting. Suppose $\gamma : H \to G$ is a map of reductive algebraic groups and $\varphi : X \to X'$ a map from an $H$-variety $X$ to a $G$-variety $X'$ such that
$\varphi(hx) = \gamma(h) \varphi(x)$. Then in $\S 6$ of \cite{bl} Bernstein and Lunts define a functor
$$
\q\varphi_* : \eqder{H}{X} \to \eqder{G}{X'}
$$
and a left adjoint
$$
\q\varphi^* :  \eqder{G}{X'} \to \eqder{H}{X}.
$$
Naturally when $\gamma$ is the identity these agree with $\varphi_*$ and $\varphi^*$. We define the equivariant intersection cohomology groups of $X$ to be 
$$
\eqih{*}{G}{X} := \shom{*+d_X}{\q\pi_*\eqic{G}{X}}
$$
where $\pi: X \to pt$ is the map to a point which is considered as a space for the trivial group, and $d_X$ the complex dimension of $X$. Note that this definition ignores the extra module structure obtained on the graded vector space $\eqih{*}{G}{X}$ from considering the point as a $G$-space.
\begin{remark}
Even for finite dimensional spaces the push-forward $\q\varphi_*$ does not preserve the bounded derived category --- this is easily seen because equivariant cohomology can be infinite dimensional --- and it is the need to use this functor which forces us to work with bounded below complexes.
\end{remark}
There are two important results which we will use. First there is an equivariant version of the famous decomposition theorem \cite[6.2.5]{bbd}. Suppose $\varphi : X \to X'$ is a proper $G$-equivariant morphism. Then we have
\begin{theorem}[{see \cite[\S 5]{bl}}]
\label{decomp thm} There is a (non-canonical) direct sum decomposition 
$$
\varphi_* \eqic{G}{X} \cong \bigoplus_{\alpha} {\imath_\alpha}_*\mathcal{IC}\udot_{\!\!\!G}(V_\alpha;\shf{L}_\alpha)[l_\alpha]
$$
where $\shf{L}_\alpha$ is an irreducible $G$-equivariant local system on the smooth part of the closed subvariety $V_\alpha$ of $X'$ and $l_\alpha \in \zz$.
\end{theorem} 
Secondly let us suppose that $1 \to N \to G \to H \to 1$ is a short exact sequence of reductive groups, that $X$ is a $G$-space upon which $N$ acts with only finite stabilisers and that further all the $N$-orbits are closed and the geometric quotient map $\varphi : X \to X/N$ is affine. Then we have
\begin{theorem}[{\cite[9.1]{bl}}]  
\label{quotient equivalence}
\begin{enumerate}
\item The functor $\q\varphi_* : \eqder{G}{X} \to \eqder{H}{X/N}$ preserves both the usual and the perverse $t$-structures so that it restricts to $\q\varphi_* : \eqsh{G}{X} \to \eqsh{H}{X/N}$ and $\q\varphi_* : \eqperv{G}{X} \to \eqperv{H}{X/N}[d_N]$ where $d_N$ is the complex dimension of $N$;
\item $\q\varphi_* \q\varphi^* = \id$;
\item $\q\varphi_* \eqic{G}{X} \cong \eqic{H}{X/N}[d_N]$.
\end{enumerate}
\end{theorem} 
In this paper we relax the conditions that $N$ acts with finite stabilisers and closed orbits and ask what should then replace the third statement above. We find that more generally $\eqic{H}{X \git N}[d_N]$ is a direct summand of $\q\varphi_*\eqic{G}{X}$.  
\section{Defining a pull-back for the quotient map}
\label{pull-back}
Suppose $L$ is an ample line bundle on a projective variety $X$ upon which a connected reductive algebraic group $G$ acts $L$-linearly. Suppose $N$ is a connected normal subgroup of $G$ and let $H$ be the quotient:
$$
1 \to N \to G \to H \to 1.
$$
We will write $X^s$ and $X^{ss}$ for the subsets of stable and semistable points with respect to the induced $N$-linearisation on $L$ (\emph{not} with respect to the $G$-linearisation). We assume that $X^s \neq \emptyset$. Consider the geometric invariant theory quotient 
$$
X^{ss} \stackrel{\varphi}{\longrightarrow} X^{ss} \git N
$$
where $X^{ss}$ is the Zariski open set of $N$-semistable points. The normality of $N$ means that $X^{ss}$ is $G$-invariant. Thus $G$ acts on $X^{ss} \git N$ via the homomorphism $G \to H$. We would like to apply the equivariant decomposition theorem to $\varphi$, which is $G$-equivariant with respect to these actions. However we cannot do so because $\varphi$ is not, in general, proper. To circumvent this difficulty we consider $\varphi$ as a rational map from $X$ to $X^{ss}\git N$ and (equivariantly) resolve the points of indeterminacy:
\[
\xymatrix{
\tilde X \ar@{-->}[drr]^{\tilde \varphi} \ar[d] & &\\
X & X^{ss} \ar@{_{(}->}[l] \ar[r]_{\varphi\quad} & X^{ss}\git N
}
\]
where $\tilde X$ is the blowup of $X$ along a closed subscheme supported on $X \setminus X^{ss}$ and $\tilde \varphi$ extends $\varphi$, see \cite[Thm. 1]{re} and cf. \cite[7.17.3]{ha}. By the equivariant decomposition theorem we then have a (non-canonical) decomposition
\begin{equation}
\label{decomp}
\tilde \varphi_* \icg{\tilde X} \cong \bigoplus_{\alpha} {\imath_\alpha}_*\mathcal{IC}\udot_{\!\!\!G}(V_\alpha;\shf{L}_\alpha)[l_\alpha].
\end{equation}
Since $\tilde \varphi$ is onto there exists a non-empty Zariski open $U \subset X^{ss} \git N$ such that upon applying the forgetful functor we have
$$
\tilde \varphi_*\ic{\tilde X} |_U \cong \bigoplus \shf{L}_\alpha|_U[l_\alpha]. 
$$
As $X^s$ is open and, by assumption, non-empty we may further assume that $U \subset \varphi(X^s)$. Since the fibres of $\varphi$ are connected $\shom{-d_X}{\varphi_*\ic{X^{ss}} |_U}$ is the constant sheaf on $U$ and the restriction
$$
\shom{-d_X}{\tilde \varphi_*\ic{\tilde X} |_U} \longrightarrow  \shom{-d_X}{\varphi_*\ic{X^{ss}} |_U}
$$
is onto. It follows that $\shom{l_\alpha-d_X}{\shf{L}_\alpha|_U}$ must also be the constant sheaf for some $\alpha$. In other words $\shf{L}_\alpha|_U$ is the constant sheaf in degree $l_\alpha-d_X$ so that
$$
{\imath_\alpha}_*\mathcal{IC}\udot_{\!\!\!G}(V_\alpha;\shf{L}_\alpha)[l_\alpha] \cong \eqic{G}{X^{ss}\git N}[d_N]
$$
is a direct summand of $\varphi_* \eqic{G}{\tilde X}$.

Let $\q_*$ be the push-forward functor $\eqder{G}{X^{ss}\git N} \to \eqder{H}{X^{ss}\git N}$. Since $N$ acts trivially on $X^{ss} \git N$  
$$
\q_*\icg{X^{ss} \git N} \cong \eqic{H}{X^{ss}\git N ; \shf{L}} 
$$
where $\shf{L}$ is a local system with stalk $H^*_N$. In particular since $H^0_N = \qq$ there is a morphism $\eqic{H}{X^{ss}\git N} \to \q_*\icg{X^{ss} \git N}$. Hence we can define a composition
\begin{equation}
\label{lambda}
\eqic{H}{X^{ss}\git N}[d_N] \to \q_*\icg{X^{ss}\git N}[d_N] \to \q \tilde \varphi_* \icg{\tilde X} \to \q \varphi_* \icg{X^{ss}}
\end{equation}
which we denote by $\lambda$. 
\begin{remark}
Intuitively we think of this as a pull-back induced by $\varphi : X^{ss} \to X^{ss} \git N$. We must be careful with this viewpoint however because $\lambda$ is not necessarily unique. The example to bear in mind is that of a variety $V$ which has two small resolutions $W_1$ and $W_2$. Let $W$ be a common resolution so that we have:
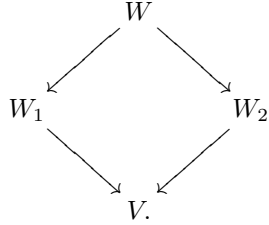
\[
\xymatrix{
& W \ar[dr] \ar[dl] &\\
W_1 \ar[dr] && W_2 \ar[dl]\\
&V.
}
\]
We know that there are natural isomorphisms $\h{*}{W_1} \cong \ih{*}{V} \cong \h{*}{W_2}$ of graded vector spaces but that the ring structures on $\h{*}{W_1}$ and $\h{*}{W_2}$ may differ so that their images in $\h{*}{W}$ cannot be the same. This shows that we cannot expect a canonical pull-back $\ih{*}{V} \to \h{*}{W}$.
\end{remark}
\begin{corollary}
\label{first cor}
The morphism $\lambda$ induces a map $\eqih{*}{H}{X^{ss}\git N} \to \eqih{*}{G}{X^{ss}}$ on hypercohomology groups.
\end{corollary}
\begin{example}
Let us suppose $X$ is a smooth projective variety and take $N=G$. There are two cases where an explicit description of a subspace of $\eqh{*}{G}{X^{ss}}$ corresponding to $\ih{*}{X^{ss}\git G}$ is already known. 
\begin{enumerate}
\item First suppose that $G=\cc^*$. Let $\mathcal{F}^0$ be the set of fixed point components of $\cc^*$ in $X^{ss}$.  For $F \in \mathcal{F}^0$ we define $N_F^+$ to be the set of $x \in X^{ss}$ such that $\lim_{t\to 0} tx \in F$, and similarly $N_F^-$ to be the set of points such that $\lim_{t\to \infty} tx \in F$. Set 
$$
c(F)=2\min(\dim_\cc N^+_F,\dim_\cc N^-_F)-1.
$$
Since $\cc^*$ acts trivially on $F$ we have $\eqh{*}{\cc^*}{F} \cong \h{*}{F} \otimes H^*_{\cc^*}$. In \cite[\S 5]{kw} it is shown that $\ih{*}{X^{ss}\git \cc^*}$ is isomorphic (as a vector space) to the set of classes in $\eqh{*}{\cc^*}{X^{ss}}$ whose restriction to $\eqh{*}{\cc^*}{F}$ has degree $< c(F)$ in the second factor of $\h{*}{F} \otimes H^{*}_{\cc^*}$ for each $F \in \mathcal{F}_0$. Furthermore it is possible to explicitly construct a morphism 
$$
\ic{X^{ss}\git \cc^*} \longrightarrow \q\varphi_* \eqc{\cc^*}{X^{ss}}
$$
which induces this identification, and to see that it is essentially unique. 
\item The second situation is when the action of $G$ (which need not now be $\cc^*$) on $X$ is \emph{weakly balanced}. This notion was introduced in \cite{kiem}. It consists of two conditions, the first of which ensures that the quotient map $\varphi : X^{ss} \to X^{ss} \git G$ is sufficiently well behaved that there is a natural pull-back $\varphi^*$ from $\ih{*}{X^{ss}\git G}$ to $\eqh{*}{G}{X^{ss}}$. In fact $\varphi$ is, in a slightly extended sense, placid --- see \cite{placid}. The second part of the weakly balanced condition allows us to identify the image of this pull-back and to show that it is in fact injective. The image is described explicitly in \cite{kiem}.
\end{enumerate}
In both the above examples a certain choice of $\lambda$ induces an inclusion of the intersection cohomology of the quotient $X^{ss} \git G$ into the $G$-equivariant cohomology of $X^{ss}$. This is not a coincidence of these examples but, as we shall see in the next section, a general feature of our situation. What perhaps is special about these examples is that there is a \emph{canonical} choice of inclusion.  
\end{example}
\section{Stable resolutions}
\label{stable}
We remind the reader that, as above, whenever we write a superscript $(s)s$ we mean (semi)stability with respect to the normal subgroup $N$ of $G$. Let us fix a choice of morphism
$$
\lambda : \eqic{H}{X^{ss}\git N}[d_N] \to \q \varphi_*\icg{X^{ss}}
$$
as above. We now show how to construct a morphism 
$$
\kappa : \q \varphi_*\icg{X^{ss}} \to \eqic{H}{X^{ss}\git N}[d_N]
$$
such that $\kappa\lambda = \id$. More precisely we show that every $N$-stable resolution of $(X,X^{ss})$ induces such a morphism.
\begin{definition}
\label{res defn}
An \emph{$N$-stable resolution} of $(X,X^{ss})$ is given by a commutative diagram of $G$-spaces
\[
\xymatrix{
Y \ar[d]_\rho & U \ar@{_(->}[l] \ar[d] \\
X & X^{ss} \ar@{_(->}[l] 
}
\]
where $\rho$ is proper and birational, and $U$ an open subset of $Y$. Further we assume that there is an ample line bundle $M \to Y$ and a $G$-linearisation on $M$ such that 
\begin{enumerate}
\item every point of $U$ is $N$-stable for $M|_U$ (but we do not necessarily assume that $U \subset Y^s$); 
\item the induced map $\sigma : U \git N \to X^{ss} \git N$ is proper and birational.
\end{enumerate}
We will usually suppress $M$ and simply write $\rho : (Y,U) \to (X,X^{ss})$ or even just $\rho$. Note that our assumption that $X^s \neq \emptyset$ implies that there is a non-empty open subset $V$ of $X^s$ such that the restriction of $\rho$ to $\rho^{-1}V \to V$ is an isomorphism. 
\end{definition}
\begin{example}
\label{res ex}
\begin{enumerate}
\item Suppose we can choose a new linearisation of the $G$ action on $X$ with the properties that $N$-stability and semistability coincide i.e. $X^{ss}_{\textrm{new}} = X^s_{\textrm{new}}$ and that $X^s_{\textrm{new}} \subset X^{ss}$. Then the inclusion 
$$(X,X^s_{\textrm{new}}) \hookrightarrow (X,X^{ss})$$ 
is an $N$-stable resolution. 
\item We reinterpret the paper \cite{lt} in the framework we have introduced. (In fact \cite{lt} is more general since it applies to symplectic reductions by $S^1$ and not just to algebraic quotients.)

Suppose $X$ is smooth and  $G=N=\cc^*$. Let $\mathcal{F}^0$ be the set of fixed point components of $\cc^*$ in $X^{ss}$. For $F \in \mathcal{F}^0$ we define $N_F^+$ to be the set of $x \in X^{ss}$ such that $\lim_{t\to 0} tx \in F$, and similarly $N_F^-$ to be the set of points such that $\lim_{t\to \infty} tx \in F$. Then
$$
X^{s} = X^{ss} \setminus \bigcup_{F \in \mathcal{F}^0} N_F^+ \cup N_F^-.
$$
We define a stable resolution by taking $Y=X$ (with $\rho = \id$) and
$$
U = X^{ss} \setminus \bigcup_{F \in \mathcal{F}^0} N_F^{d(F)}
$$
where $d(F)=\pm$ depending upon whether $\dim_\cc N^+_F$ or $\dim_\cc N^-_F$ is the larger, if they are equal then we make an arbitrary choice. (Note that in this case $U \not \subset X^s$ but that every point of $U$ is nevertheless stable for the restriction of $L$ to $U$.) The induced map $U \git \cc^* \to X^{ss}\git \cc^*$ is shown to be a \emph{small} resolution in \cite{lt}. This follows from an application of the ideas of \cite{h1} to a $G$-invariant neighbourhood of each $F\in \mathcal{F}^0$. Hence $\ih{*}{X^{ss}\git \cc^*} \cong \ih{*}{U \git \cc^*} \cong \eqh{*}{\cc^*}{U}$. Lerman and Tolman then go on to compute the kernel of the restriction $\eqh{*}{\cc^*}{X^{ss}} \to \eqh{*}{\cc^*}{U}$, and so express $\ih{*}{X^{ss}\git \cc^*}$ as a quotient of $\eqh{*}{\cc^*}{X^{ss}}$. 
\item 
\label{kirwans res}
Suppose $X$ is smooth and take $N=G$. In \cite{k2} Kirwan describes a canonical partial desingularisation of the quotient $X^{ss}\git G$. This is constructed by taking the quotient of a projective variety $Y$ (with a suitably linearised $G$ action) obtained from $X$ by a sequence of blowups. Initially we blow up along the closure in $X$ of the smooth $G$-invariant subvariety
$$
GZ^{ss}_R  = \{ x \in X^{ss} \ |\  \Stab_G x \textrm{ conjugate to } R\}
$$ 
where $R$ is a connected reductive subgroup of $G$ such that $GZ_R^{ss}$ has maximal codimension amongst all such subvarieties. There is an induced action of $G$ on the blowup $\pi : \tilde X \to X$ which can be linearised on the ample bundle formed by twisting the pull-back of a sufficiently large power of $L$ with minus the exceptional divisor. It is shown in \cite[Lemma 6.1]{k2} that $\pi(\tilde X^{ss})\subset X^{ss}$ and $G\tilde Z_R^{ss} = \emptyset$. Continuing inductively we can construct a $G$-stable resolution.

We can refine this procedure slightly by working relative to a nontrivial normal subgroup $N$. Now we consider blowing up along the closures of subvarieties of the form $GZ_R^{ss}$ where $ss$ refers to $N$-semistability and $R$ is a connected reductive subgroup of $N$. The same procedure will then construct an $N$-stable resolution. 
\end{enumerate}
\end{example}
\begin{proposition}
\label{existence of stable}
There is always at least one $N$-stable resolution. 
\end{proposition}
\begin{proof}
To construct a resolution  we first of all $G$-equivariantly resolve the singularities of $X$ (see \cite[Thm. 1]{re}). Thus we have a smooth $G$-variety $\tilde X$ and a $G$-equivariant map $\tilde X \to X$ which factors as a finite sequence of blowups of smooth $G$-invariant subvarieties. 

Suppose the first of these blowups is
$\pi : \hat X \to X$. The action of $G$ on $\hat X$ can be linearised on an ample line bundle of the form $\pi^* L^d \otimes \mathcal{O}(-E)$ where $E$ is the exceptional divisor and $d$ is sufficiently large. As usual let $\hat X^{ss}$ denote the $N$-semistable points of $\hat X$ with respect to this linearisation. I claim that $\pi(\hat X^{ss}) \subset X^{ss}$. To see this note that if $x \in \hat X^{ss}$ then there is a $N$-invariant section $\sigma = \sigma_1 \otimes \sigma_2$ of $\pi^* L^d \otimes \mathcal{O}(-E)$ with $\sigma(x) \not = 0$. Since $\mathcal{O}(-E)$ is trivial away from $E$ we deduce that $\sigma_1$ is $N$-invariant and, of course, $\sigma_1(x) \not = 0$. Equivalently $\pi(x)$ is $N$-semistable for the obvious linearisation on $L^d$. But by the final remark of \cite[Chapter1, \S 5]{mfk} the $N$-semistable points for the linearisations on $L$ and $L^d$ coincide. 

Proceeding inductively we see that we can linearise the action of $G$ on $\tilde X$ in such a way that we have a map $(\tilde X, \tilde X^{ss}) \to (X,X^{ss})$. We can now apply Kirwan's resolution, relative to $N$, to the smooth $G$-variety $\tilde X$ (see example \ref{kirwans res} of \ref{res ex}). In this way we will obtain an $N$-stable resolution of $(X,X^{ss})$.
\end{proof}
Let $\rho$ be an $N$-stable resolution of $X^{ss}$ so that, in the notation of \ref{res defn},  we have a diagram:
\[
\xymatrix{
Y \ar_{\rho}[d] & U \ar[d] \ar@{_(->}[l] \ar[r]^{\psi} & U \git N \ar[d]^\sigma\\
X & X^{ss} \ar@{_(->}[l] \ar[r]_\varphi & X^{ss}\git N .
}
\]
We can apply the equivariant decomposition theorem to the proper maps $\rho$ and $\sigma$ to obtain morphisms
$$
\icg{X^{ss}} \longrightarrow  \rho_* \icg{U} 
\textrm{\quad and \quad}
\sigma_* \eqic{H}{U \git N} \longrightarrow \eqic{H}{X^{ss} \git N}.
$$
These are of course not canonical but their restrictions to $V$ and $\varphi (V)$ respectively are the natural quasi-isomorphisms induced by $\rho$ and $\rho^{-1}$. By theorem \ref{quotient equivalence} $\q\psi_*\icg{U} \cong \eqic{H}{U\git N}[d_N]$ because $N$ acts with only finite stabilisers on $U$. So we can compose these morphisms to obtain
\begin{equation}
\label{defn of kappa}
\kappa: \q\varphi_* \icg{X^{ss}} \longrightarrow \eqic{H}{X^{ss}\git N}[d_N].
\end{equation}
\begin{theorem}
\label{main thm}
The composition  $\kappa\lambda[-d_N]$ is the identity on $\eqic{H}{X^{ss}\git N}$.
\end{theorem}
\begin{proof}
We know from \cite[\S 5]{bl} that $\eqic{H}{X^{ss}\git N}$ is a simple object in the heart of the perverse $t$-structure on $\eqder{H}{X^{ss}\git N}$ so that 
$$
\Hom{\eqic{H}{X^{ss}\git N}}{\eqic{H}{X^{ss}\git N}} \cong \qq.
$$
We can easily check that $\kappa\lambda[-d_N]$ restricts to the identity on $\varphi(V)$ and hence must be the identity. 
\end{proof}
\begin{corollary}
\label{main cor}
The map $\eqih{*}{H}{X^{ss}\git N} \hookrightarrow \eqih{*}{G}{X^{ss}}$ induced by $\lambda$ is an inclusion. Any $N$-stable resolution can be used to define (not necessarily uniquely) a projection $\eqih{*}{G}{X^{ss}}\to \eqih{*}{H}{X^{ss}\git N}$ which is split by this inclusion.  
\end{corollary}
\begin{proof}
This follows immediately from proposition \ref{existence of stable} and theorem \ref{main thm}.
\end{proof}
\begin{remark}
In \cite{k3} the partial desingularisation constructed in \cite{k2} (see example \ref{res ex} part \ref{kirwans res}) is used to define, as we have above, a map $\eqh{*}{G}{X^{ss}} \to \ih{*}{X^{ss}\git G}$ when $X$ is smooth. \cite[Thm. 2.5]{k3} states that this map is a surjection. This is then used to give an algorithm for computing the intersection Betti numbers of $X^{ss} \git G$. 

Unfortunately the given argument that the map is surjective is incorrect due to a confusion on page 499 where it says ``by 2.24 there is a commutative diagram
\[
\xymatrix{
H^*_G(E)/H^*_N(Z^{ss}_R) \ar[r] \ar[d] & H^*_G(Y) \ar[d]^{\hat{\phi}_Y}\\
I\!H^*(E \git G)/I\!H^*(\mathcal{N}\git G) \ar[r] & I\!H^*(Y \git G)
}
\]
where the horizontal maps are the inclusions associated to these decompositions.'' In fact Kirwan's 2.24 shows that the same four objects form a commutative diagram but where the horizontal maps arise from Gysin sequences and are not those associated to ``these decompositions''.

Corollary \ref{main cor} provides a new proof that Kirwan's map is surjective and hence re-validates the algorithm.  
\end{remark}

\bibliography{full}

\begin{thebibliography}{MFK94}

\bibitem[BBD82]{bbd}
A.~Beilinson, J.~Bernstein, and P.~Deligne.
\newblock Faisceaux pervers.
\newblock {\em Ast\'erisque}, 100, 1982.
\newblock Proc. C.I.R.M. conf\'erence: Analyse et topologie sur les espaces
  singuliers.

\bibitem[BL94]{bl}
J.~Bernstein and V.~Lunts.
\newblock {\em Equivariant sheaves and functors}.
\newblock Number 1578 in Lecture Notes in Mathematics. Springer--Verlag, 1994.

\bibitem[GK96]{gk}
V.~Guillemin and J.~Kalkman.
\newblock The {J}effrey--{K}irwan localization theorem and residue operations
  in equivariant cohomology.
\newblock {\em Journal fur die Reine und Angewandte Mathematik}, 470:123--142,
  1996.

\bibitem[GM83]{gm2}
M.~Goresky and R.~Macpherson.
\newblock Intersection homology theory {II}.
\newblock {\em Inventiones Mathematicae}, 71:77--129, 1983.

\bibitem[G{\'{o}}m99]{gom}
T.~G{\'{o}}mez.
\newblock Algebraic stacks.
\newblock arXiv:math.AG/9911199, 1999.

\bibitem[Har77]{ha}
R.~Hartshorne.
\newblock {\em Algebraic Geometry}.
\newblock Springer--Verlag, 1977.

\bibitem[Hu92]{h1}
Y.~Hu.
\newblock The geometry and topology of quotient varieties of torus actions.
\newblock {\em Duke {M}ath. {J}.}, 68(1):151--184, 1992.

\bibitem[JK95]{jk}
L.~Jeffrey and F.~Kirwan.
\newblock Localization for non-abelian groups.
\newblock {\em Topology}, 34(2):291--327, 1995.

\bibitem[Kie98]{kiem}
Y-H. Kiem.
\newblock Intersection cohomology of quotients of nonsingular varieties.
\newblock Preprint, 1998.

\bibitem[Kie00]{placid}
Y.-H. Kiem.
\newblock Note on placid maps.
\newblock Preprint, October 2000.

\bibitem[Kir84]{k1}
F.~Kirwan.
\newblock {\em Cohomology of Quotients in Symplectic and Algebraic Geometry}.
\newblock Number~34 in Mathematical Notes. Princeton University Press, 1984.

\bibitem[Kir85]{k2}
F.~Kirwan.
\newblock Partial desingularisations of quotients of nonsingular varieties and
  their {B}etti numbers.
\newblock {\em Annals of Mathematics}, 122:41--85, 1985.

\bibitem[Kir86]{k3}
F.~Kirwan.
\newblock Rational intersection homology of quotient varieties.
\newblock {\em Inventiones Mathematicae}, 86:471--505, 1986.

\bibitem[KW00]{kw}
Y-H. Kiem and J.~Woolf.
\newblock The cosupport axiom, equivariant cohomology and the intersection
  cohomology of certain symplectic quotients.
\newblock Preprint, 2000.

\bibitem[LT00]{lt}
E.~Lerman and S.~Tolman.
\newblock Intersection cohomology of ${S}^1$ symplectic quotients and small
  resolutions.
\newblock {\em Duke Math J.}, 103(1):79--99, 2000.

\bibitem[MFK94]{mfk}
D.~Mumford, J.~Fogarty, and F.~Kirwan.
\newblock {\em Geometric invariant theory}.
\newblock Springer--Verlag, third edition, 1994.

\bibitem[RY99]{re}
Z.~Reichstein and B.~Youssin.
\newblock Equivariant resolution of points of indeterminacy.
\newblock arXiv:math.AG/0006099, 1999.

\bibitem[TW98]{tw}
S.~Tolman and J.~Weitsman.
\newblock The cohomology rings of abelian symplectic quotients.
\newblock arXiv:math.DG/9812006, 1998.

\bibitem[Wit92]{wi}
E.~Witten.
\newblock Two dimensional gauge theories revisited.
\newblock {\em J. Geom. Phys.}, 9:303--368, 1992.

\end{thebibliography}
\bibliographystyle{alpha}

\end{document}